\input amstex
\documentstyle{amsppt}
\magnification=\magstep1
 \hsize 13cm \vsize 18.35cm \pageno=1
\loadbold \loadmsam
    \loadmsbm
    \UseAMSsymbols
\topmatter
\NoRunningHeads
\title A Note on $q$-Bernstein polynomials
 \endtitle
\author
  Taekyun Kim \endauthor
 \keywords $q$-Bernstein polynomials, $q$-Euler numbers , $q$-Stirling numbers, fermionic $p$-adic integrals
\endkeywords

\abstract
Recently, Simsek-Acikgoz([17]) and Kim-Jang-Yi([9]) have studied  the $q$-extension of Bernstein polynomials. In this paper
we propose the $q$-extension of Bernstein polynomials of degree $n$, which are different $q$-Bernstein polynomials of Simsek-Acikgoz([17])
and Kim-Jang-Yi([9]).  From these $q$-Bernstein polynomials, we derive some fermionic $p$-adic integral representations of  several $q$-Bernstein
 type polynomials. Finally, we investigate some identities between $q$-Bernstein polynomials and $q$-Euler numbers.
\endabstract
\thanks  2000 AMS Subject Classification: 11B68, 11S80, 60C05, 05A30
\newline  The present Research has been conducted by the research
Grant of Kwangwoon University in 2010
\endthanks
\endtopmatter

\document

{\bf\centerline {\S 1. Introduction}}

 \vskip 20pt
Let $C[0, 1]$ denote the set of continuous function on $[0, 1]$. For $f\in C[0, 1]$, Bernstein introduced the following  well known linear operators
(see [1, 3]):
$$ \Bbb B_n(f|x)=\sum_{k=0}^n f(\frac{k}{n})\binom{n}{k}x^k(1-x)^{n-k}=
\sum_{k=0}^n f(\frac{k}{n})B_{k,n}(x). \tag1 $$
Here $\Bbb B_n(f|x)$ is called Bernstein operator of order $n$ for $f$. For $k, n \in \Bbb Z_{+}(=\Bbb N \cup \{0\})$, the Bernstein polynomials of degree $n$ is defined by
$$B_{k,n}(x)=\binom{n}{k}x^k(1-x)^{n-k}, \text{ (see [1, 2, 3]).}\tag2$$
A Bernoulli trial involves performing an experiment once and noting whether
 a particular event $A$ occurs. The outcome of Bernoulli trial is said to be ``success"  if $A$ occurs and a ``failure" otherwise.
Let $k$ be the number of successes in $n$ independent Bernoulli trials, the probabilities of $k$ are given by the binomial probability law:
$$p_n(k)=\binom{n}{k}p^k(1-p)^{n-k},  \text{ for $k=0, 1, \cdots, n$,}$$
where $p_n(k)$ is the probability of $k$ successes in $n$ trials. For example, a communication system transmit binary information over channel
that introduces random bit errors with probability $\xi=10^{-3}$. The transmitter transmits each information bit three times, an a decoder
takes a majority vote of the received bits to decide on what the transmitted bit was. The receiver can correct a single error, but it will make
the wrong decision if the channel introduces two or more errors. If we view each transmission as a Bernoulli trial in which a ``success" corresponds to the introduction of an error, then the probability of two or more errors in three Bernoulli trial is
$$p (k\geq 2)=\binom{3}{2}(0.001)^2(0.999)+\binom{3}{3}(0.001)^3\approx 3(10^{-6}), \text{ see [18]}.$$
 By the definition of Bernstein polynomials(see Eq.(1) and Eq.(2)), we can see that
Bernstein basis is the probability mass function of binomial distribution. In the reference [15] and [16], Phillips proposed a generalization
of classical Bernstein polynomials based on $q$-integers. In the last decade some new generalizations of well known positive linear operators based on $q$-integers were introduced and studied by several authors(see [1-21]). Let $0<q<1$. Define the $q$-numbers of $x$ by
$ [x]_q =\frac {1-q^x}{1-q}$(see [1-21]). Recently, Simsek-Acikgoz([17]) and Kim-Jang-Yi([9]) have studied the $q$-extension of Bernstein polynomials, which are different Phillips $q$-Bernstein polynomials.  Let $p$ be a fixed odd prime number. Throughout this paper $\Bbb Z_p$, $\Bbb Q_p$, and $\Bbb C_p$ denote the rings of $p$-adic integers, the fields of $p$-adic rational numbers, and the completion of algebraic closure of $\Bbb Q_p$, respectively. The $p$-adic absolute value in $\Bbb C_p$ is normalized in such way that $|p|_p=\frac{1}{p}.$ As well known definition, Euler polynomials are defined by
$$ \frac{2}{e^t+1}e^{xt}=\sum_{n=0}^{\infty}E_n(x)\frac{t^n}{n!}, \text{ (see [1-14]).}\tag3$$
In the special case, $x=0$, $E_n(0)=E_n$ are called the $n$-th Euler numbers. By (3), we see that the recurrence formula of Euler numbers is given by
$$E_0=1, \text{ and }  (E+1)^n+E_n=0  \text{  if  $n>0$, (see [12]),}\tag4$$
with the usual convention of replacing  $E^n$ by $E_n$.
When one talks of $q$-analogue, $q$ is variously considered as an indterminate, a complex number $q\in \Bbb C$, or a $p$-adic
number $q\in \Bbb C_p$. If $q\in \Bbb C$, we normally assume $|q|<1$. If $q\in \Bbb C_p$, we normally always assume that $|1-q|_p<1$.
As the $q$-extension of (4), author defined the $q$-Euler numbers as follows:
$$ E_{0,q}=1, \text{ and } (qE_q+1)^n+E_{n,q}=0 \text{ if $ n>0,$ ( see [21]),} \tag5$$
 with the usual convention of replacing $E_q^n$ by $E_{n,q}$.
Let $UD(\Bbb Z_p)$ be the space of uniformly differentiable function on $\Bbb Z_p$. For $f\in UD(\Bbb Z_p)$, the fermionic $p$-adic $q$-integral was defined by
$$ I_q(f)=\int_{\Bbb Z_p}f(x)d\mu_{-q}(x)=\lim_{N\rightarrow \infty}
\frac{1}{1+q^{p^N}}\sum_{x=0}^{p^N-1} f(x)(-q)^x, \text{ (see [12])}.\tag6$$
In the special case, $q=1$,  $I_1(f)$ is called the fermionic $p$-adic integral on $\Bbb Z_p$ (see [12, 21]). By (6) and the definition of $I_1(f)$, we see that
$$ I_1(f_1)+ I_1(f)=2f(0), \text{ where $f_1(x)=f(x+1).$}\tag7$$
For $n\in \Bbb N$, let $f_n(x)=f(x+n)$. Then we can also see that
$$ I_1(f_n)+(-1)^{n-1}I_1(f)=2\sum_{l=0}^{n-1}(-1)^{n-l-1}f(l), \text{ (see [21])}.\tag8$$
From (5), (7) and (8), we note that
$$\int_{\Bbb Z_p}e^{[x]_qt}d\mu_{-1}(x)=\sum_{n=0}^{\infty}E_{n,q}\frac{t^n}{n!}
=\sum_{n=0}^{\infty}\left( \frac{2}{(1-q)^n}\sum_{l=0}^n\binom{n}{l}\frac{(-1)^l}{1+q^l}\right)\frac{t^n}{n!}.\tag9$$
Thus we have
$$ E_{n,q}= \frac{2}{(1-q)^n}\sum_{l=0}^n\binom{n}{l}\frac{(-1)^l}{1+q^l},  \text{ (see [21])}.$$
In [21], the $q$-Euler polynomials are defined by
$$E_{n,q}(x)=\int_{\Bbb Z_p}[y+x]_q^n d\mu_{-1}(x)=\frac{1}{(1-q)^n}\sum_{l=0}^n\binom{n}{l}(-1)^l \frac{q^{lx}}{
1+q^l}.\tag10$$
By (9) and (100, we get
$$E_{n,q}(x)=\sum_{l=0}^n\binom{n}{l}q^{lx}E_{l,q}=(q^xE_q+1)^n,\tag11$$
with the usual convention of replacing $E_q^n$ by $E_{n,q} $.
In this paper we firstly consider the $q$-Bernstein polynomials of degree $n$ in $\Bbb R$, which are different
$q$-Bernstein polynomials of Simsek-Acikgoz([17]) and Kim-Jang-Yi([9]).
From these $q$-Bernstein polynomials, we try to study for the fermionic $p$-adic integral representations of
the several $q$-Bernstein type polynomials on $\Bbb Z_p$. Finally, we give some interesting identities between $q$-Bernstein
polynomials and $q$-Euler numbers.

\vskip 20pt

{\bf\centerline {\S 2. $q$-Bernstein  Polynomials}} \vskip 10pt
   For $n, k \in \Bbb Z_{+},$ the generating function for $B_{k,n}(x)$ is introduced by Acikgoz and Araci as follows:
   $$F^{(k)}(t,x)=\frac{te^{(1-x)t}x^k}{k!}=\sum_{n=0}^{\infty}B_{k,n}(x)\frac{t^n}{n!}, \text{ (see [1, 9, 10, 17])}.\tag12$$
   For  $k, n\in\Bbb Z_{+}$,  $0<q<1$ and $x\in [0,  1]$,
   consider the $q$-extension of (12) as follows:
   $$\aligned
   F_q^{(k)}(t,x)&=\frac{(t[x]_q)^ke^{[1-x]_{\frac{1}{q}}t}}{k!}=\frac{[x]_q^k}{k!}\sum_{n=0}^{\infty}\frac{[1-x]_{\frac{1}{q}}^n}{n!}
   t^{n+k}=\sum_{n=k}^{\infty}\left(\frac{n![x]_q^k[1-x]_{\frac{1}{q}}^{n-k}}{(n-k)!k!}\right)\frac{t^n}{n!}\\
   &= \sum_{n=k}^{\infty}\binom{n}{k}[x]_q^k [1-x]_{\frac{1}{q}}^{n-k}\frac{t^n}{n!}=\sum_{n=k}^{\infty}B_{k,n}(x,q)\frac{t^n}{n!} .
   \endaligned\tag13$$
Because $B_{k, 0}(x, q)=B_{k,1}(x,q)=\cdots=B_{k, k-1}(x,q)=0$, we obtain the following generating function for $B_{k, n}(x,q)$:
   $$F_q^{(k)}(t,x)=\frac{(t[x]_q)^ke^{[1-x]_{\frac{1}{q}}t}}{k!}=\sum_{n=0}^{\infty}B_{k,n}(x,q)\frac{t^n}{n!}, \text{  where $k\in\Bbb Z_{+}$
   and $x\in [0, 1]$.}$$
Thus, for $ k, n\in \Bbb Z_{+}$, we note that
$$\aligned
 B_{k, n}(x, q) &=\binom{n}{k}[x]_q^k [1-x]_{\frac{1}{q}}^{n-k},  \text{ if $ n\geq k$ },\\
 &= 0,  \text{ if $k<n$.}
\endaligned\tag14$$
By (14), we easily get $\lim_{q\rightarrow 1}B_{k, n}(x,q)=B_{k, n}(x).$ For $0\leq k \leq n$, we have
$$\aligned
&[1-x]_{\frac{1}{q}}B_{k, n-1}(x, q)+[x]_qB_{k-1, n-1}(x, q)\\
&=[1-x]_{\frac{1}{q}}\binom{n-1}{k}[x]_q^k[1-x]_{\frac{1}{q}}^{n-k-1}
+[x]_q\binom{n-1}{k-1}[x]_q^{k-1}[1-x]_{\frac{1}{q}}^{n-k}\\
&=\binom{n-1}{k}[x]_q^k[1-x]_{\frac{1}{q}}^{n-k}+\binom{n-1}{k-1}[x]_q^k[1-x]_{\frac{1}{q}}^{n-k}
=\binom{n}{k}[x]_q^k[1-x]_{\frac{1}{q}}^{n-k},
\endaligned$$
and the derivative of the $q$-Bernstein polynomials of degree $n$ are also polynomials of degree $n-1$.
$$\aligned
&\frac{d}{dx}B_{k-1,n}(x, q)\\
&=k\binom{n}{k}[x]_q^{k-1}[1-x]_{\frac{1}{q}}^{n-k}\left(\frac{\log q}{q-1}\right)q^x +
\binom{n}{k}[x]_q^k(n-k)[1-x]_{\frac{1}{q}}^{n-k-1}\left(\frac{\log q}{1-q}\right)q^{x}\\
&= \frac{\log q}{q-1}q^x\left( n\binom{n-1}{k-1}[x]_q^{k-1}[1-x]_{\frac{1}{q}}^{n-k}-n\binom{n-1}{k}[x]_q^k[1-x]_{\frac{1}{q}}^{n-1-k} \right)\\
&=n \left( B_{k-1, n-1}(x, q)-B_{k, n-1}(x, q)\right)\frac{\log q}{q-1} q^x.
\endaligned$$

Therefore,  we obtain the following theorem.

\proclaim{ Theorem 1} For $k, n\in \Bbb Z_{+}$ and $x\in [0, 1]$, we have
$$[1-x]_{\frac{1}{q}}B_{k, n-1}(x, q)+ [x]_qB_{k-1, n-1}(x, q)=B_{k, n}(x, q), $$
and
$$\frac{d}{dx}B_{k,n}(x,q)=n\left(B_{k-1, n-1}(x, q)-B_{k, n-1}(x, q) \right)\frac{\log  q}{ q-1}q^x.$$
\endproclaim
Let $f$ be a continuous function on $[0, 1]$. Then the $q$-Bernstein operator of order $n$ for $f$ is defined by
$$\Bbb B_{n,q}(f|x)=\sum_{k=0}^n f(\frac{k}{n})B_{k, n}(x, q), \text{ where $0\leq x \leq 1 $ and $ n \in \Bbb Z_{+}$}.
 \tag15$$
By (14) and (15), we see that
$$\Bbb B_{n, q}(1|x)=\sum_{k=0}^nB_{k, n}(x, q)=\sum_{k=0}^n\binom{n}{k}[x]_q^k[1-x]_{\frac{1}{q}}^{n-k}=\left([x]_q+[1-x]_{\frac{1}{q}}\right)^n=1.$$
Also, we get from (15) that for $f(x)=x,$
$$\Bbb B_{n, q}(x|x)=\sum_{k=0}^n\frac{k}{n}\binom{n}{k}[x]_q^k[1-x]_{\frac{1}{q}}^{n-k}=\sum_{k=0}^{n-1}\binom{n-1}{k}[x]_q^{k+1}[1-x]_{\frac{1}{q}}^{n-k-1}
=[x]_q .$$
The $q$-Bernstein polynomials are symmetric polynomials in the following sense:
$$B_{n-k,n}(1-x, \frac{1}{q})=\binom{n}{n-k}[1-x]_{\frac{1}{q}}^{n-k}[x]_q^k=B_{k,n}(x, q). $$
Thus, we obtain the following theorem.

\proclaim{ Theorem 2}
For $n, k \in \Bbb Z_{+}$ and $(x\in [0, 1]$, we have
$$B_{n-k, n}(1-x, \frac{1}{q})=B_{k, n}(x, q).$$
Moreover,
$\Bbb B_{n, q}(1|x)=1$ and   $\Bbb B_{n, q}(x|x)=[x]_q.$
\endproclaim
  From (15), we note that
  $$\aligned
  \Bbb B_{n,q}(f|x)&=\sum_{k=0}^n f(\frac{k}{n})B_{k, n}(x, q)=\sum_{k=0}^nf(\frac{k}{n})\binom{n}{k}[x]_q^k[1-x]_{\frac{1}{q}}^{n-k}\\
  &=\sum_{k=0}^n f(\frac{k}{n})\binom{n}{k}[x]_q^k\sum_{j=0}^{n-k}\binom{n-k}{j}(-1)^j[x]_q^j.
  \endaligned\tag16$$
   By the definition of binomial coefficient,  we easily get
   $$\binom{n}{k}\binom{n-k}{j}=\binom{n}{k+j}\binom{k+j}{k}.$$
 Let $k+j=m$. Then we have
 $$\binom{n}{k}\binom{n-k}{j}=\binom{n}{m}\binom{m}{k}. \tag17$$
  From (16) and (17), we have
  $$\Bbb B_{n, q}(f|x)=\sum_{m=0}^n\binom{n}{m}[x]_q^m\sum_{k=0}^m\binom{m}{k}(-1)^{m-k}f(\frac{k}{n}). \tag18$$
 Therefore, we obtain the following theorem.
 \proclaim{ Theorem 3}
 For $f\in C[0, 1]$ and $n\in\Bbb Z_{+}$, we have
 $$\Bbb B_{n, q}(f|x)=\sum_{m=0}^n\binom{n}{m}[x]_q^m\sum_{k=0}^m\binom{m}{k}(-1)^{m-k}f(\frac{k}{n}). $$
 \endproclaim
It is well known that the second kind  stirling numbers are defined by
$$\frac{(e^t-1)^k}{k!}=\frac{1}{k!}\sum_{l=0}^k \binom{k}{l}(-1)^{k-l}e^{lt}
=\sum_{n=0}^{\infty}s(n,k)\frac{t^n}{n!}, \text{ for $k\in\Bbb N$, (see [12, 21]).} \tag 19$$
Let $\Delta $ be the shift difference operator with $\Delta f(x)=f(x+1)-f(x).$ By iterative process, we easily get
$$\Delta^n f(0)=\sum_{k=0}^n\binom{n}{k}(-1)^{n-k}f(k). \tag20$$
From (19) and (20), we can easily derive the following equation (21).
$$\frac{1}{k!}\Delta^k0^n=s(n, k). \tag21$$
By (18) and (20) we obtain the following  theorem.
\proclaim{ Theorem 4}
For $f\in C[0, 1]$ and $n\in \Bbb Z_{+}$, we have
$$\Bbb B_{n, q}(f|x)=\sum_{k=0}^n \binom{n}{k}[x]_q^k \Delta^k f(\frac{0}{n}). $$
\endproclaim
In the special case, $f(x)=x^m$($m\in\Bbb Z_{+}$), we have the following corollary.

\proclaim{ Corollary 5}
For $x\in [0, 1]$ and $m, n\in \Bbb Z_{+}$, we have
$$n^m \Bbb B_{n, q}(x^m|x)=\sum_{k=0}^n \binom{n}{k}[x]_q^k\Delta^k0^m, $$
and
$$n^m\Bbb B_{n, q}(x^m|x)=\sum_{k=0}^n \binom{n}{k}[x]_q^kk!s(m, k).$$
\endproclaim

For $x, t \in \Bbb C$ and $n\in\Bbb Z_{+}$ with $n \geq k$, consider
$$\frac{n!}{2\pi i}\int_{C}\frac{([x]_q t)^k}{k!}e^{([1-x]_{\frac{1}{q}}t)}\frac{dt}{ t^{n+1}}, \tag22$$
where $C$ is a circle around the origin and integration is in the positive  direction.
We see from the definition of the $q$-Bernstein  polynomials and the basic theory of complex analysis including Laurent series that
$$ \int_C\frac{([x]_q t)^k}{k!}e^{([1-x]_{\frac{1}{q}} t)}\frac{dt}{t^{n+1}}
=\sum_{m=0}^{\infty}\int_C \frac{B_{k, m}(x, q)t^m}{m!} \frac{dt}{t^{n+1}}=2\pi i\left(\frac{B_{k, n}(x, q)}{n!}\right).\tag23 $$
We get from (22) and (23) that
$$\frac{n!}{2\pi i}\int_C \frac{([x]_q t)^k}{k!}e^{([1-x]_{\frac{1}{q}} t)}\frac{dt}{t^{n+1}} = B_{k, n}(x, q),\tag24 $$
and
$$\aligned
\int_C\frac{([x]_q t)^k}{k!}e^{([1-x]_{\frac{1}{q}} t)}\frac{dt}{t^{n+1}}&=\frac{[x]_q^k}{k!}\sum_{m=0}^{\infty}
\left(\frac{[1-x]_{\frac{1}{q}}^m}{m!}\int_C t^{m-n-1+k}dt \right)\\
&=2\pi i
\left( \frac{[x]_q^k[1-x]_{\frac{1}{q}}^{n-k}}{k!(n-k)!}\right)=\frac{2\pi i}{n!}\binom{n}{k}[x]_q^k[1-x]_{\frac{1}{q}}^{n-k}.\\
\endaligned \tag25$$
By (22) and (25), we see that
$$ \frac{n!}{2\pi i}\int_C\frac{([x]_q t)^k}{k!}e^{([1-x]_{\frac{1}{q}} t)}\frac{dt}{t^{n+1}}=\binom{n}{k}[x]_q^k[1-x]_{\frac{1}{q}}^{n-k}.
\tag26$$
From (24) and (26), we note that
$$ B_{k, n}(x, q)=\binom{n}{k}[x]_q^k[1-x]_{\frac{1}{q}}^{n-k}.$$
By the definition of $q$-Bernstein polynomials, we easily get
$$\aligned
&\left(\frac{n-k}{n}\right)B_{k, n}(x,q)+\left(\frac{k+1}{n}\right)B_{k+1, n}(x, q) \\
&=\left(\frac{ (n-1)!}{k!(n-k-1)!}\right)[x]_q^k[1-x]_{\frac{1}{q}}^{n-k}
+\left(\frac{(n-1)!}{k!(n-k-1)!}\right)[x]_q^{k+1}[1-x]_{\frac{1}{q}}^{n-k-1}\\
&=\left([1-x]_{\frac{1}{q}}+[x]_q \right)B_{k, n-1}(x, q)=B_{k, n-1}(x, q).
\endaligned$$
Therefore, we can write $q$-Bernstein  polynomials as a linear combination of polynomials of higher order.
\proclaim{ Theorem 6}
For $k, n \in \Bbb Z_{+}$ and $x\in [0, 1]$, we have
$$ \left(\frac{n+1-k}{n+1}\right)B_{k, n+1}(x,q)+\left(\frac{k+1}{n+1}\right)B_{k+1, n+1}(x, q)=B_{k, n}(x, q).$$
\endproclaim
We easily get from (14) that for $n, k \in \Bbb N,$
$$\aligned
&\left( \frac{n-k+1}{k}\right)\left(\frac{[x]_q}{[1-x]_{\frac{1}{q}}}\right)B_{k-1, n}(x, q)\\
&=\left( \frac{n-k+1}{k}\right)\left(\frac{[x]_q}{[1-x]_{\frac{1}{q}}}\right)\binom{n}{k-1}[x]_q^{k-1}[1-x]_{\frac{1}{q}}^{n-k+1}\\
&=\left(\frac{n!}{k!(n-k)!}\right)[x]_q^k [1-x]_{\frac{1}{q}}^{n-k}=B_{k, n}(x, q).
\endaligned$$
Therefore, we obtain the following corollary.
\proclaim{ Corollary 7}
For $k, n\in \Bbb N$ and $x\in [0, 1]$, we have
$$ \left( \frac{n-k+1}{k}\right)\left(\frac{[x]_q}{[1-x]_{\frac{1}{q}}}\right)B_{k-1, n}(x, q)=B_{k, n}(x, q).$$
\endproclaim
By (14) and binomial theorem, we easily see that
$$B_{k, n}(x,q)=\binom{n}{k}[x]_q^k\sum_{l=0}^{n-k}\binom{n-k}{l}(-1)^l[x]_q^l=\sum_{l=k}^n\binom{l}{k}\binom{n}{k}(-1)^{l-k}[x]_q^l .$$
Therefore, we obtain the following theorem.
\proclaim{ Theorem 8}
For $k, n \in \Bbb Z_{+}$ and $x\in [0, 1]$, we have
$$B_{k, n}(x, q)=\sum_{l=k}^n\binom{l}{k}\binom{n}{k}(-1)^{l-k}[x]_q^l .$$
\endproclaim
It is possible to write $[x]_q^k$ as a linear combination of the $q$-Bernstein polynomials by using the degree  evaluation formulae
and mathematical induction.  We easily see from the property of the $q$-Bernstein polynomials that
$$\sum_{k=1}^n \left(\frac{k}{n}\right)B_{k, n}(x,q)=\sum_{k=0}^{n-1}\binom{n-1}{k}[x]_q^{k+1}[1-x]_{\frac{1}{q}}^{n-k-1}
=[x]_q, $$
and that
$$\sum_{k=2}^n\frac{\binom{k}{2}}{\binom{n}{2}}B_{k, n}(x, q)=\sum_{k=0}^{n-2}\binom{n-2}{k}[x]_q^{k+2}[1-x]_{\frac{1}{q}}^{n-2-k}
=[x]_q^2. $$
Continuing this process, we get
$$\sum_{k=j}^n \frac{\binom{k}{j}}{\binom{n}{j} }B_{k, n}(x, q)=[x]_q^j, \text{ for $j\in \Bbb Z_{+}$}.$$
Therefore, we obtain the following theorem.
\proclaim{ Theorem 9}
For $n, j \in \Bbb Z_{+}$ and $ x\in [0, 1]$, we have
$$\sum_{k=j}^n \frac{\binom{k}{j}}{\binom{n}{j} }B_{k, n}(x, q)=[x]_q^j.$$
\endproclaim
In [7], the $q$-stirling numbers of the second kind are defined by
$$s_q(n, k)=\frac{ q^{-\binom{k}{2}}}{[k]_q!}\sum_{j=0}^k (-1)^j q^{\binom{j}{2}}{\binom{k}{j}}_q[k-j]_q^n, \tag27$$
where ${\binom{k}{j}}_q=\frac{[k]_q!}{[j]_q![k-j]_q!}$ and $[k]_q!=\prod_{i=1}^k [i]_q$. For $n\in \Bbb Z_{+}$, it is known that
$$[x]_q^n=\sum_{k=0}^n q^{\binom{k}{2}}{\binom{x}{k}}_q[k]_q!s_q(n, k), \text{ (see [7, 21]).}\tag28$$
By (27), (28) and Theorem 7, we obtain the following corollary.
\proclaim{ Corollary 10}
For $n, j \in \Bbb Z_{+}$ and $x\in [0, 1]$, we have
$$\sum_{k=j}^n \frac{\binom{k}{j}}{\binom{n}{j} }B_{k, n}(x, q) =\sum_{k=0}^j q^{\binom{k}{2}}{\binom{x}{k}}_q[k]_q!s_q(j, k). $$
\endproclaim
\vskip 20pt
{\bf\centerline {\S 3. On fermionic $p$-adic integral representations of $q$-Bernstein polynomials}} \vskip 10pt

In this section we assume that $q\in \Bbb C_p$ with $|1-q|_p<1$. From (10) we note that
$$ E_{n, \frac{1}{q}}(1-x)=\int_{\Bbb Z_p}[1-x+x_1]_{\frac{1}{q}}^n d\mu_{-1}(x_1)
=(-1)^nq^n\int_{\Bbb Z_p}[x+x_1]_q^n d\mu_{-1}(x_1), \text{ (see [21])}.\tag29$$
From (29) we have
$$\aligned
&\int_{\Bbb Z_{p}}[1-x]_{\frac{1}{q}}^n d\mu_{-1}(x)=q^n(-1)^n\int_{\Bbb Z_p}[x-1]_q^n d\mu_{-1}(x)\\
&=\int_{\Bbb Z_p}\left( 1-[x]_q\right)^n d\mu_{-1}(x)=(-1)^nq^n E_{n, \frac{1}{q}}(-1)=E_{n, q}(2).
\endaligned$$
By (5) and (10), we easily get
$$E_{n, q}(2)=2+E_{n, q}, \text{ if  $n >0$.} $$
Thus,  we  obtain the following theorem.
\proclaim{ Theorem 11}
For $n \in \Bbb N$, we have
$$\int_{\Bbb Z_p}[1-x]_{\frac{1}{q}}^nd\mu_{-1}(x)=\int_{\Bbb Z_p}\left(1-[x]_q \right)^nd\mu_{-1}(x)=2+\int_{\Bbb Z_p}[x]_q^n d\mu_{-1}(x).$$
\endproclaim
By using Theorem 11, we derive our main results in this section. Taking the fermionic $p$-adic integral on $\Bbb Z_p$ for one $q$-Bernstein polynomials in (14), we get
$$\aligned \int_{\Bbb Z_p}B_{k, n}(x,q)d\mu_{-1}(x)&=\binom{n}{k}\int_{\Bbb Z_p}[x]_q^k[1-x]_{\frac{1}{q}}^{n-k} d\mu_{-1}(x)\\
&=\binom{n}{k}\sum_{l=0}^{n-k}\binom{n-k}{l}(-1)^l\int_{\Bbb Z_p}[x]_q^{k+l}d\mu_{-1}(x)\\
&=\binom{n}{k}\sum_{l=0}^{n-k}\binom{n-k}{l}(-1)^l E_{k+l, q}.
\endaligned\tag30$$
From (14) and Theorem 2, we note that
$$\aligned
\int_{\Bbb Z_p}B_{k,n}(x, q)d\mu_{-1}(x)&=\int_{\Bbb Z_p}B_{n-k, n}(1-x, \frac{1}{q}) d\mu_{-1}(x)\\
&=\binom{n}{k}\sum_{j=0}^k\binom{k}{j}(-1)^{k+j}\int_{\Bbb Z_p}[1-x]_{\frac{1}{q}}^{n-j} d\mu_{-1}(x).
\endaligned\tag31$$
For $n>k$, by (31) and Theorem 11, we get
$$\aligned
\int_{\Bbb Z_p} B_{k, n}(x, q) d\mu_{-1}(x)&
=\binom{n}{k}\sum_{j=0}^k \binom{k}{j}(-1)^{k+j}\left(2+\int_{\Bbb Z_p}[x]_q^{n-j}d\mu_{-1}(x) \right)\\
&= 2+ E_{n, q}, \text{ if $k=0$ }\\
&=\binom{n}{k}\sum_{j=0}^k \binom{k}{j}(-1)^{k+j}E_{n-j, q}, \text{ if $k>0$}.
\endaligned\tag32$$
From $m, n, k \in \Bbb Z_{+}$ with $m+n>2k$, the fermionic $p$-adic integral for multiplication of two $q$-Bernstein polynomials on $\Bbb Z_{p}$
can be given by the following relation:
$$\aligned
&\int_{\Bbb Z_p} B_{k, n}(x,q) B_{k, m}(x,q) d\mu_{-1}(x)
=\binom{n}{k}\binom{m}{k}\int_{\Bbb Z_{p}}[x]_q^{2k}[1-x]_{\frac{1}{q}}^{n+m-2k} d\mu_{-1}(x)\\
&=\binom{n}{k}\binom{m}{k}\sum_{j=0}^{2k}\binom{2k}{j}(-1)^{j+2k}\int_{\Bbb Z_p}[1-x]_{\frac{1}{q}}^{n+m-j}d\mu_{-1}(x)\\
&=\binom{n}{k}\binom{m}{k}\sum_{j=0}^{2k}(-1)^{j+2k}\left( 2+\int_{\Bbb Z_p}[x]_q^{n+m-j} d\mu_{-1}(x)\right).
\endaligned\tag33$$
From (33), we have
$$\aligned
\int_{\Bbb Z_p} B_{k, n}(x,q) B_{k, m}(x,q) d\mu_{-1}(x)&= 2+E_{n+m, q}, \text{  if $k=0$}\\
&=\binom{n}{k}\binom{m}{k}\sum_{j=0}^{2k}\binom{2k}{j}(-1)^{j+2k}E_{n+m-j, q}, \text{ if $k>0$}.
\endaligned$$
For $m, k \in \Bbb Z_{+}$, it is difficult to show that
$$\int_{\Bbb Z_p} B_{k, n}(x,q) B_{k, m}(x,q) d\mu_{-1}(x)=\binom{n}{k}\binom{m}{k}
\sum_{j=0}^{n+m-2k}\binom{n+m-2k}{j}(-1)^j E_{j+2k, q}. \tag34$$
Continuing this process we obtain the following theorem.
\proclaim{ Theorem 12}
(I). For $n_1, \cdots, n_s, k \in \Bbb Z_{+}$ ($s\in \Bbb N$) with $n_1+\cdots+n_s>sk$, we have
$$\int_{\Bbb Z_p}\left(\prod_{i=1}^s B_{k, n_i}(x, q) \right) d\mu_{-1}(x)=2+E_{n_1+\cdots+n_s, q}, \text{ if $k=0$},$$
and
$$\int_{\Bbb Z_p}\left(\prod_{i=1}^s B_{k, n_i}(x, q) \right) d\mu_{-1}(x)=\prod_{i=1}^s\binom{n_i}{k}\sum_{j=0}^{sk}\binom{sk}{j}(-1)^{sk-j}E_{n_1+\cdots+n_s-j,q}, \text{ if $k>0$}.$$
(II). Let $k, n_1, \cdots n_s\in \Bbb Z_{+}$ ($s\in \Bbb N $). Then we have
$$\aligned
&\int_{\Bbb Z_{p}}\left(\prod_{i=1}^s B_{k, n_i}(x, q)\right)d\mu_{-1}(x)\\
& =\left(\prod_{i=1}^s \binom{n_i}{k}\right)
\sum_{j=0}^{\sum_{i=1}^sn_i-sk} \binom{\sum_{i=1}^sn_i-sk}{j} (-1)^j E_{j+sk, q}.\endaligned $$
\endproclaim
By Theorem 12, we obtain the following corollary.
\proclaim{ Corollary 13} For $n_1, \cdots, n_s, k \in \Bbb Z_{+}$ ($s\in \Bbb N$) with $n_1+\cdots+n_s>sk$, we have
$$\sum_{j=0}^{\sum_{i=1}^sn_i-sk} \binom{\sum_{i=1}^sn_i-sk}{j} (-1)^j E_{j+sk, q}= 2+E_{n_1+\cdots+n_s, q}, \text{ if $k=0$}, $$
and
$$\aligned
&\sum_{j=0}^{\sum_{i=1}^sn_i-sk} \binom{\sum_{i=1}^sn_i-sk}{j} (-1)^j E_{j+sk, q}\\
&=\sum_{j=0}^{sk}\binom{sk}{j}(-1)^{sk-j}E_{n_1+\cdots+n_s-j,q}, \text{ if $k>0$}.\endaligned$$
\endproclaim
Let $m_1, \cdots, m_s, n_1, \cdots, n_s, k \in\Bbb Z_{+}$ ($s\in \Bbb N$) with $m_1n_1+\cdots +m_sn_s> (m_1+\cdots+m_s)k$.
By the definition of $B_{k, n_s}^{m_s}(x, q)$, we can also easily see that
$$\aligned
&\int_{\Bbb Z_p}\left( \prod_{i=1}^s B_{k, n_i}^{m_i}(x, q)\right)d\mu_{-1}(x)\\
&=\prod_{i=1}^s {\binom{n_i}{k}}^{m_i}\sum_{j=0}^{k\sum_{i=1}^s m_i}\binom{k\sum_{i=1}^s m_i}{j}
(-1)^{k\sum_{i=1}^s m_i-j}
\int_{\Bbb Z_p}[1-x]_{\frac{1}{q}}^{\sum_{i=1}^s n_im_i-j} d\mu_{-1}(x)\\
&=\prod_{i=1}^s {\binom{n_i}{k}}^{m_i}\sum_{j=0}^{k\sum_{i=1}^s m_i}\binom{k\sum_{i=1}^s m_i}{j}
(-1)^{k\sum_{i=1}^s m_i-j}
\left(2+E_{\sum_{i=1}^sm_i n_i-j, q} \right).
\endaligned$$

\vskip 20pt

\quad Taekyun Kim

\quad Division of General Education-Mathematics, Kwangwoon
University, Seoul

\quad 139-701, S. Korea
 e-mail:\text{ tkkim$\@$kw.ac.kr}

\enddocument